\documentclass[12pt]{article}
\usepackage{amsmath}
\usepackage{epsf}
\usepackage{amssymb}
\usepackage{graphicx}
\usepackage{ifpdf}
\usepackage{caption}

\usepackage{array}
\usepackage{booktabs} %调整表格线与上下内容的间隔
\usepackage{multirow}

\newtheorem{theorem}{Theorem}[section]
\newtheorem{claim}{Claim}[section]
\newtheorem{definition}{Definition}[section]
\newtheorem{lemma}{Lemma}[section]

\begin{document}
\pagestyle{empty}
\renewcommand{\thefootnote}{\fnsymbol{footnote}}

\begin{center}
{\bf \Large The generalized 3-connectivity of burnt pancake graphs and godan graphs}\footnote{This
work was supported by National Science Foundation of China(No.12271157), Natural Science Foundation of Hunan Province (No.2022JJ30028) and Hunan Education Department Foundation(No.21C0762).}
\vskip 5mm

{{\bf Jing Wang$^1$, Zuozheng Zhang$^1$, Yuanqiu Huang$^2$}\\[2mm]
$^1$ School of Mathematics, Changsha University, Changsha 410022, China\\
$^2$ School of Mathematics, Hunan Normal University, Changsha 410081, China}\\[6mm]
\end{center}
\date{}

\noindent{\bf Abstract}\; The generalized $k$-connectivity of a graph $G$, denoted by $\kappa_k(G)$, is the minimum number of internally edge disjoint $S$-trees for any $S\subseteq V(G)$ and $|S|=k$. The generalized $k$-connectivity is a natural extension of the classical connectivity and plays a key role in applications related to the modern interconnection networks. The burnt pancake graph $BP_n$ and the godan graph $EA_n$ are two kinds of Cayley graphs which posses many desirable properties. In this paper, we investigate the generalized 3-connectivity of $BP_n$ and $EA_n$. We show that $\kappa_3(BP_n)=n-1$ and $\kappa_3(EA_n)=n-1$.

\noindent{\bf Keywords} interconnection network, Cayley graph, generalized $k$-connectivity, tree \\
{\bf MR(2000) Subject Classification} 05C40, 05C05

\section{Introduction}
\label{secintro}

With rapid development and advances of very large scale integration technology and wafer-scale integration technology, multiprocessor systems have been widely designed and used in our daily life. It is well known that the underlying topology of the multiprocessor systems can be modelled by a connected graph $G=(V(G),E(G))$, where $V(G)$ is the set of processors and $E(G)$ is the set of communication links of multiprocessor systems.

As Cayley graphs have a lot of properties which are desirable in an interconnection network, such as vertex transitivity, edge transitivity, hierarchical structure, high fault tolerance, many researchers have proposed Cayley graphs as models for interconnection networks \cite{Heydemann1997}. Let $\mathcal{A}$ be a finite group and $\Omega$ be a subset of $\mathcal{A}$, where the identity of the group does not belong to $\Omega$. The {\it Cayley graph} ${\rm Cay}(\mathcal{A}, \Omega)$ is a digraph with vertex set $\mathcal{A}$ and arc set $\{(a, a\cdot \omega)| a\in \mathcal{A}, \omega\in \Omega \}$. If $\Omega=\Omega^{-1}$, then ${\rm Cay}(\mathcal{A}, \Omega)$ is an undirected graph, where $\Omega^{-1}=\{\omega^{-1}| \omega\in \Omega \}$.

Fault tolerance has become increasingly significant nowadays since multiprocessor systems failure is inevitable. The connectivity  is a key parameter for  measuring fault tolerance of the network. A subset $S\subseteq V(G)$ of a connected graph $G$ is called a {\it vertex-cut} if $G\setminus S$ is disconnected or trivial. The {\it connectivity} $\kappa(G)$ of $G$ is defined as the minimum cardinality over all vertex-cuts of $G$. Note that the larger $\kappa(G)$ is, the more reliable the network is. A well known theorem of Whitney \cite{Whitney1932} provides an equivalent definition of connectivity. For each 2-subset $S=\{x,y\}\subseteq V(G)$, let $\kappa(S)$ denote the maximum number of internally disjoint ($x,y$)-paths in $G$. Then
\begin{equation*}
\kappa(G)=\min\{\kappa(S) |S\subseteq V(G)\; {\rm and} \; |S|=2\}.
\end{equation*}

The generalized $k$-connectivity, which was introduced by Chartrand et al. \cite{Chartrand1984}, is a strengthening of connectivity and can be served as an essential parameter for measuring reliability and fault tolerance of the network. Let $G=(V(G),E(G))$ be a simple graph, $S$ be a subset of $V(G)$. A tree $T$ in $G$ is called an $S$-{\it tree}, if $S\subseteq V(T)$. The trees $T_1, T_2, \cdots, T_r$ are called {\it internally edge disjoint $S$-trees} if $V(T_i)\cap V(T_j)=S$ and $E(T_i)\cap E(T_j)=\emptyset$ for any integers $1\le i\ne j\le r$. $\kappa_G(S)$ denote the maximum number of internally edge disjoint $S$-trees. For an integer $k$ with $2\le k\le |V(G)|$, the {\it generalized $k$-connectivity} of $G$, denoted by $\kappa_k(G)$, is defined as
\begin{equation*}
\kappa_k(G)=\min\{\kappa_G(S) |S\subseteq V(G)\; {\rm and} \; |S|=k\}.
\end{equation*}

The generalized 2-connectivity is exactly the classical connectivity. Over the past few years, research on the generalized connectivity has received meaningful progress. Li et al. \cite{SLi2012n} derived that it is NP-complete for a general graph $G$ to decide whether there are $l$ internally disjoint trees connecting $S$, where $l$ is a fixed integer and $S\subseteq V(G)$. Authors in \cite{HZLi2014,SLi2010} investigated the upper and lower bounds of the generalized connectivity of a general graph $G$.

Many authors tried to study exact values of the generalized connectivity of graphs. The generalized $k$-connectivity of the complete graph, $\kappa_k(K_n)$, was determined in \cite{Chartrand2010} for every pair $k,n$ of integers with $2\le k\le n$. The generalized $k$-connectivity of the complete bipartite graphs $K_{a,b}$ are obtained in \cite{SLi2012b} for all $2\le k\le a+b$. Apart from these two results, the generalized $k$-connectivity of other important classes of graphs, such as, Cartesian product graphs \cite{HZLi2012,HZLi2017}, hypercubes \cite{HZLi2012,SLin2017}, several variations of hypercubes \cite{ZhaoHao20191,ZhaoHao20192,Wei2021,Wang2021}, Cayley graphs \cite{SLi2017,ZHao20193,SLi2016,Hao20191,ANn12018,ANn12019,ANn2018}, have draw many scholars' attention. So far, as we can see, the results on the generalized $k$-connectivity of networks are almost about $k=3$.

The burnt pancake graph $BP_n$\cite{Chin2009} and the godan graph $EA_n$ \cite{GodanIn2022,godangraph2022} are two kinds of Cayley graphs, each of which is regular, vertex-transitive and used to design various commercial multiprocessor machines. In this paper, we try to evaluate the reliability and fault tolerance of the former two networks by studying their generalized 3-connectivity. Section \ref{secpreli} introduces some notations and definitions. The generalized 3-connectivity of the burnt pancake graph $BP_n$ and the $n$-dimensional godan graphs $EA_n$ are investigated in Section \ref{secBPn} and Section \ref{secEAn}, respectively.

\section{Preliminaries}\label{secpreli}

This section is dedicated to introduce some necessary preliminaries. We only consider a simple, connected graph $G=(V(G),E(G))$ with $V(G)$ be its vertex set and $E(G)$ be its edge set. For a vertex $x\in V(G)$, the {\it degree} of $x$ in $G$, denoted by ${\rm deg}_G(x)$, is the number of edges of $G$ incident with $x$. Denote $\delta(G)$ the {\it minimum degree} of vertices of $G$. We can abbreviate $\delta(G)$ to $\delta$ if there is no confusion. A graph is $d$-{\it regular} if ${\rm deg}_G(x)=d$ for every vertex $x\in V(G)$. For a vertex $x\in V(G)$, we use $N_G(x)$ to denote the neighbour vertices set of $x$ and $N_G[x]$ to denote $N_G(x)\cup \{x\}$. Let $V'\subseteq V(G)$, denote by $G\backslash V'$ the graph obtained from $G$ by deleting all the vertices in $V'$ together with their incident edges. Denote by $G[V']$ the subgraph of $G$ induced on $V'$.

Let $P$ be a path in $G$ with $x$ and $y$ be its two terminal vertices, then $P$ is called an $(x,y)$-{\it path}. Two ($x,y$)-paths $P_1$ and $P_2$ are {\it internally disjoint} if they have no common internal vertices, that is, $V(P_1)\cap V(P_2)=\{x, y\}$.

Li et al. \cite{SLi2010} gave upper and lower bounds of $\kappa_3(G)$ for a general graph $G$.

\begin{lemma}\label{lemupperK}(\cite{SLi2010})
Let $G$ be a connected graph with minimum degree $\delta$. If there are two adjacent vertices of degree $\delta$, then $\kappa_3(G)\le \delta-1$.
\end{lemma}

\begin{lemma}\label{lemkr}(\cite{SLi2010})
Let $G$ be a connected graph with $n$ vertices. For every two integers $k$ and $r$ with $k\ge 0$ and $r\in\{0,1,2,3\}$, if $\kappa(G)=4k+r$, then $\kappa_3(G)\ge 3k+\lceil \frac{r}{2}\rceil$.
\end{lemma}

\begin{lemma}\label{lemxypath} (\cite{Bondy})
Let $G$ be a $k$-connected graph, and let $x$ and $y$ be a pair of distinct vertices of $G$. Then there exists $k$ internally disjoint ($x,y$)-paths in $G$.
\end{lemma}

\begin{lemma}\label{lemKfan} (\cite{Bondy})
Let $G$ be a $k$-connected graph, let $x$ be a vertex of $G$ and let $Y\subseteq V(G)-\{x\}$ be a set of at least $k$ vertices of $G$. Then there exists a $k$-fan in $G$ from $x$ to $Y$, that is, there exists a family of $k$ internally disjoint ($x,Y$)-paths whose terminal vertices are distinct in $Y$.
\end{lemma}

\begin{lemma}\label{lemXYpaths} (\cite{Bondy})
Let $G$ be a $k$-connected graph, and let $X$ and $Y$ be subsets of $V(G)$ of cardinality at least $k$. Then there exists a family of $k$ pairwise disjoint ($X,Y$)-paths in $G$.
\end{lemma}

\section{The generalized 3-connectivity of the burnt pancake graph}\label{secBPn}

Let $[n]=\{1,2,\cdots, n\}$.  Let $S_n$ be a symmetric group on $[n]$. For an integer $i$, it is well known that $|i|$ denotes the absolute value of $i$. Denote $\bar{i}=-i$ in this paper. Let $[[n]]$ be the set $[n]\cup \{\bar{i} ~|~ i\in[n]\}$. A {\it signed permutation} of $[n]$ is an $n$-permutation $x_1x_2\cdots x_n$ of $[[n]]$ such that $|x_1||x_2|\cdots |x_n|$ forms a permutation of $[n]$. For a signed permutation $x=x_1x_2\cdots x_n$ of $[[n]]$ and an integer $i$ ($1\le i\le n$), the $i$th {\it signed prefix reversal} of $x$ is denoted by $x^i=\bar{x}_i \bar{x}_{i-1}\cdots \bar{x}_1x_{i+1}\cdots x_n$.

\begin{definition}\label{defBPn}(\cite{Chin2009})
For $n\ge 2$, an $n$-dimensional burnt pancake graph $BP_n$ is a graph with vertex set $V(BP_n)=\{x\, | x$ is a signed permutation of $[[n]]\}$. Two vertices $x=x_1x_2\cdots x_n$ and $y=y_1y_2\cdots y_n$ are adjacent in $BP_n$ if and only if there exists an integer $i$ ($1\le i\le n$) such that $x^i=y$.
\end{definition}

\begin{figure}[htbp]
\begin{minipage}[t]{0.3\linewidth}
\centering
\resizebox{0.7\textwidth}{!} {\includegraphics{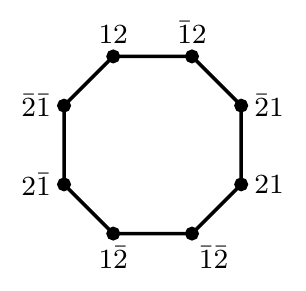}}
\caption{\small The burnt pancake graph $BP_2$} \label{figBP2}
\end{minipage}
%\quad
\begin{minipage}[t]{0.7\linewidth}
\centering
\resizebox{0.9\textwidth}{!} {\includegraphics{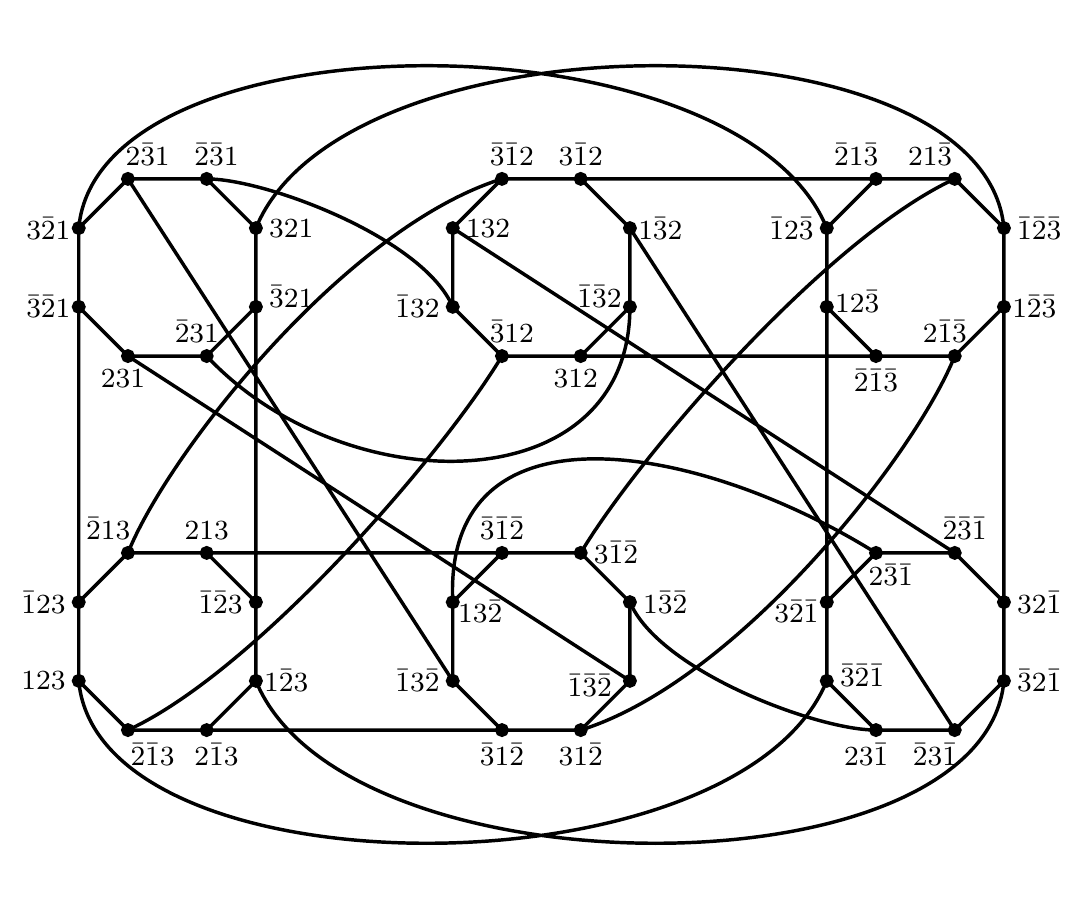}}
\caption{\small The burnt pancake graph $BP_3$} \label{figBP3}
\end{minipage}
\end{figure}

The burnt pancake graph  $BP_2$ and $BP_3$ are depicted in Figure \ref{figBP2} and Figure \ref{figBP3}, respectively. Let $Z_2\wr S_n$ be the wreath product of the cyclic group $Z_2=\{0,1\}$ and the symmetric group $S_n$ and $\Omega=\{(\underbrace{1,1,\cdots ,1,}_i0,\cdots ,0; i(i-1)\cdots 1(i+1)\cdots n)|2\le i \le n\}$. Song et al. proved in \cite{Song2015} that $BP_n$ is isomorphic to the Cayley graph Cay$(Z_2\wr S_n,~\Omega)$.

Note that, by fixing the symbol $i$ in the rightmost position of each vertex for $i\in[[n]]$, $BP_n$ can be decomposed into $2n$ vertex disjoint subgraphs $BP_n^i$, called {\it clusters}. Obviously, $BP_n^i$ is isomorphic to $BP_{n-1}$ for each $i\in[[n]]$. We write the construction of $BP_n$ symbolically as
$$BP_n=BP_n^1\oplus BP_n^{\bar{1}}\oplus BP_n^2 \oplus BP_n^{\bar{2}}\oplus \cdots \oplus  BP_n^n \oplus BP_n^{\bar{n}}.$$
An edge $(x,y)$ of $BP_n$ is called a {\it cross edge} if $x$ and $y$ belong to different clusters, moreover, $y$ is called the {\it out-neighbour} of $x$ and is denoted by $\hat{x}$ under this circumstance. For $\{i,j\}\subseteq [[n]]$, we denote by $E(BP_n^i,BP_n^j)$ the set of cross edges between $BP_n^i$ and $BP_n^j$.

\begin{lemma}\label{lemBPn1}(\cite{Song2015,Gate1979})
For $n\ge 2$, the burnt pancake graph $BP_n$ has the following properties:\\
(1) it is an $n$-regular Cayley graph with $2^nn!$ vertices and $n\times n! \times 2^{n-1}$ edges;\\
(2) $\kappa(BP_n)=n$;\\
(3) For $\{i,j\}\subseteq [[n]]$, it has $|E(BP_n^i,BP_n^j)|=(n-2)!\times 2^{n-2}$ if $i\ne \bar{j}$, otherwise, $|E(BP_n^i,BP_n^j)|=0$.
\end{lemma}

\begin{lemma}\label{lemBPn2}
For $n\ge 2$, let $H=BP_n\setminus V(BP_n^j)$, where $BP_n^j$ is any cluster of $BP_n$ for $j\in [[n]]$. Then $\kappa(H)=n-1$.
\end{lemma}

\noindent {\bf Proof}\; W.l.o.g., assume that $H=BP_n\setminus V(BP_n^{\bar{n}})$. Firstly, it is clear that $\kappa(H)\le\delta(H)=n-1$. To obtain the reverse inequality, we only need to prove that there exist ($n-1$) internally disjoint ($x,y$)-paths for any two vertices $x$ and $y$ in $H$.

If $x$ and $y$ belong to a same cluster of $H$, say $BP_n^1$. Then there are ($n-1$) internally disjoint ($x,y$)-paths in $BP_n^1$ since $\kappa(BP_n^1)=\kappa(BP_{n-1})=n-1$. Now consider that $x$ and $y$ belong to different clusters, say $BP_n^l$ and $BP_n^k$.

\vskip 2mm
{\bf Case 1.}\; $l\ne \bar{k}$.

By Lemma \ref{lemBPn1}, it has $|E(BP_n^l,BP_n^k)|=(n-2)!\times 2^{n-2}$. Choose ($n-1$) different vertices $u_1, \cdots, u_{n-1}$ from $BP_n^l$ such that their out-neighbours $\hat{u}_1, \cdots, \hat{u}_{n-1}$ belong to $BP_n^k$. This can be done since $(n-2)!\times 2^{n-2}\ge n-1$ for $n\ge 2.$

Let $U=\{u_1, \cdots, u_{n-1}\}$ and $\hat{U}=\{\hat{u}_1, \cdots, \hat{u}_{n-1}\}$. It is clear that $|U|=|\hat{U}|=n-1$. For $1\le i\le n-1$, there exists an ($n-1$)-fan $L_1, \cdots, L_{n-1}$ in $BP_n^l$ from $x$ to $U$ such that $u_i\in V(L_i)$ and an ($n-1$)-fan $Q_1, \cdots, Q_{n-1}$ in $BP_n^k$ from $y$ to $\hat{U}$ such that $\hat{u}_i\in V(Q_i)$.

For $1\le i\le n-1$, let $P_i=L_i\cup Q_i\cup \{u_i\hat{u}_i\}$. Then ($n-1$) internally disjoint ($x,y$)-paths are constructed.

\vskip 2mm
{\bf Case 2.}\; $l=\bar{k}$.

W.l.o.g., assume that $l=1$ and $k=\bar{1}$. According to Lemma \ref{lemBPn1},
$$|E(BP_n^1,BP_n^2)| =|E(BP_n^{\bar{1}},BP_n^2)|=(n-2)!\times 2^{n-2}.$$
Therefore, we can choose ($n-1$) different vertices $u_1, \cdots, u_{n-1}$ from $BP_n^1$ such that their out-neighbours $\hat{u}_1, \cdots, \hat{u}_{n-1}$ belong to $BP_n^2$ and ($n-1$) different vertices $\omega_1, \cdots, \omega_{n-1}$ from $BP_n^{\bar{1}}$ such that their out-neighbours $\hat{\omega}_1, \cdots, \hat{\omega}_{n-1}$ belong to $BP_n^2$ since $(n-2)!\times 2^{n-2}\ge n-1$ for $n\ge 2$.

Let $\hat{U}=\{\hat{u}_1, \cdots, \hat{u}_{n-1}\}$ and $\hat{W}=\{\hat{\omega}_1, \cdots, \hat{\omega}_{n-1}\}$. It is known that $\hat{U}\cap \hat{W}=\emptyset$ and $|\hat{U}|= |\hat{W}|=n-1$. According to Lemma \ref{lemXYpaths}, there are ($n-1$) disjoint ($\hat{U},\hat{W}$)-paths $R_1, \cdots, R_{n-1}$ in $BP_n^2$ such that $\{\hat{u}_i, \hat{\omega}_i\}\subseteq V(R_i)$. Likewise, there is an ($n-1$)-fan $L_1, \cdots, L_{n-1}$ from $x$ to $\{u_1, \cdots, u_{n-1}\}$ in $BP_n^1$ such that $u_i\in V(L_i)$ and an ($n-1$)-fan $Q_1, \cdots, Q_{n-1}$ from $y$ to $\{\omega_1, \cdots, \omega_{n-1}\}$ in $BP_n^{\bar{1}}$ such that $\omega_i\in V(Q_i)$ ($1\le i\le n-1$).

For $1\le i\le n-1$, let $P_i=L_i\cup Q_i\cup R_i\cup\{u_i\hat{u}_i,\omega_i\hat{\omega}_i\}$. Then ($n-1$) internally disjoint ($x,y$)-paths are obtained. The proof is done. \hfill$\Box$

\begin{lemma}\label{lemBPn3}
Let $BP_n=BP_n^1\oplus BP_n^{\bar{1}}\oplus BP_n^2 \oplus BP_n^{\bar{2}}\oplus \cdots \oplus  BP_n^n \oplus BP_n^{\bar{n}}$ for $n\ge 2$. Let $x=x_1\cdots x_{n-1}j$ be any vertex of $BP_n^{j}$ ($j\in [[n]]$). Then the out-neighbours of vertices in $N_{BP_n^{j}}[x]$ belong to $n$ different clusters of $BP_n$.
\end{lemma}

\noindent{\bf Proof}\; By Definition \ref{defBPn}, we have $N_{BP_n^{j}}[x]=\{x^i=\bar{x}_i \bar{x}_{i-1}\cdots \bar{x}_1x_{i+1}\cdots x_{n-1}j| 1\le i\le n-1\}\cup \{x\}$. Thus, the out-neighbour of $x^i$ is $\bar{j}\bar{x}_{n-1} \cdots\bar{x}_{i+1} x_1\cdots x_{i-1}x_i$, which belongs to $V(BP_n^{x_i})$, $1\le i\le n-1$. Moreover, $\hat{x}=\bar{j}\bar{x}_{n-1} \cdots\bar{x}_1\in V(BP_n^{\bar{x}_1})$. Obviously, the out-neighbours of $N_{BP_n^{j}}[x]$ belong to $n$ different clusters of $BP_n$.  \hfill$\Box$

\begin{lemma}\label{lemBPn4}
Let $BP_n=BP_n^1\oplus BP_n^{\bar{1}}\oplus BP_n^2 \oplus BP_n^{\bar{2}}\oplus \cdots \oplus  BP_n^n \oplus BP_n^{\bar{n}}$ for $n\ge 3$. Let $x,y$ and $z$ be any three vertices that belong to three different clusters of $BP_n$. Then there exist ($n-1$)-internally edge disjoint $\{x,y,z\}$-trees in $BP_n$.
\end{lemma}

\noindent{\bf Proof}\; W.l.o.g., assume that $x\in V(BP_n^{j_1})$, $y\in V(BP_n^{j_2})$ and $z\in V(BP_n^{j_3})$, where $\{j_1,j_2,j_3\}\subseteq [[n]]$. The following three cases are discussed.

\vskip 2mm
{\bf Case 1.}\; $n\ge 5$.

Since $2n-6\ge n-1$ for $n\ge 5$, there exist ($n-1$) different integers $l_1, \cdots, l_{n-1}$ such that $\{j_1,j_2,j_3,\bar{j}_1,\bar{j}_2,\bar{j}_3\}\cap\{l_1, \cdots, l_{n-1}\}=\emptyset$. By Lemma \ref{lemBPn1}, $|E(BP_n^a,BP_n^b)|=(n-2)!\times 2^{n-2}$, where $a\in\{j_1,j_2,j_3\}$ and $b\in \{l_1, \cdots, l_{n-1}\}$. For $1\le i\le n-1$, we may choose vertices $u_i\in V(BP_n^{j_1})$, $r_i\in V(BP_n^{j_2})$ and $\omega_i\in V(BP_n^{j_3})$ such that their out-neighbours $\hat{u}_i,\hat{r}_i,\hat{\omega}_i$ belong to $V(BP_n^{l_i})$. There is an
$\{\hat{u}_i,\hat{r}_i,\hat{\omega}_i\}$-tree $\widetilde{T}_i$ in $BP_n^{l_i}$ since $BP_n^{l_i}$ is connected.

Let $U=\{u_1,\cdots, u_{n-1}\}$, $R=\{r_1,\cdots, r_{n-1}\}$ and $W=\{\omega_1,\cdots, \omega_{n-1}\}$. Obviously, $|U|=|R|=|W|=n-1$. By Lemma \ref{lemKfan}, there is an ($n-1$)-fan $P_1, \cdots, P_{n-1}$ from $x$ to $U$ in $BP_n^{j_1}$ such that $u_i\in V(P_i)$ for $1\le i\le n-1$. Likewise, for $1\le i\le n-1$, there is an ($n-1$)-fan $Q_1, \cdots, Q_{n-1}$ from $y$ to $R$ in $BP_n^{j_2}$ such that $r_i\in V(Q_i)$ and an ($n-1$)-fan $L_1, \cdots, L_{n-1}$ from $z$ to $W$ in $BP_n^{j_3}$ such that $\omega_i\in V(L_i)$.

Let $T_i=\widetilde{T}_i\cup P_i\cup Q_i\cup L_i\cup \{u_i\hat{u}_i, r_i\hat{r}_i, \omega_i\hat{\omega}_i\}$, where $1\le i\le n-1$. Then $T_1, \cdots, T_{n-1}$ are ($n-1$)-internally edge disjoint $\{x,y,z\}$-trees in $BP_n$.

\vskip 2mm
{\bf Case 2.}\; $n=4$.

\vskip 2mm
{\bf Subcase 2.1.}\; There exist two integers in $\{j_1, j_2, j_3\}$, say $j_1$ and $j_2$, such that $j_1+j_2=0$.

W.l.o.g., assume that $x\in V(BP_n^1)$, $y\in V(BP_n^{\bar{1}})$ and $z\in V(BP_n^2)$.

By Lemma \ref{lemBPn1}, there exist vertex subsets $\{u_1,u_2,u_3\}\subseteq V(BP_n^1)$, $\{r_1,r_2,r_3\}\subseteq V(BP_n^{\bar{1}})$ and $\{w_1,w_2,w_3\}\subseteq V(BP_n^2)$, for which $\{\hat{u}_1, \hat{r}_1, \hat{w}_1\}\subseteq V(BP_n^3)$, $\{\hat{u}_2, \hat{r}_2, \hat{w}_2\}\subseteq V(BP_n^{\bar{3}})$ and $\{\hat{u}_3, \hat{r}_3, \hat{w}_3\}\subseteq V(BP_n^4)$. Since $BP_n^3$ is connected, there is an $\{\hat{u}_1, \hat{r}_1, \hat{w}_1\}$-tree $\widetilde{T}_1$ in $BP_n^3$. Likewise, there is an $\{\hat{u}_2, \hat{r}_2, \hat{w}_2\}$-tree $\widetilde{T}_2$ in $BP_n^{\bar{3}}$ and an $\{\hat{u}_3, \hat{r}_3, \hat{w}_3\}$-tree $\widetilde{T}_3$ in $BP_n^4$, respectively.

Moreover, there is a 3-fan $P_1, P_2, P_3$ from $x$ to $\{u_1, u_2, u_3\}$ in $BP_n^1$ for which $u_i\in V(P_i)$, $1\le i\le 3$. Analogously, for $1\le i\le 3$, there is a 3-fan $R_1,R_2,R_3$ from $y$ to $\{r_1,r_2,r_3\}$ in $BP_n^{\bar{1}}$ such that $r_i\in V(R_i)$ and a 3-fan $L_1,L_2,L_3$ from $z$ to $\{w_1,w_2,w_3\}$ in $BP_n^{2}$ such that $w_i\in V(L_i)$.

For $1\le i\le 3$, let
$$T_i=\widetilde{T}_i\cup P_i\cup R_i\cup L_i\cup \{u_i\hat{u}_i, r_i\hat{r}_i, w_i\hat{w}_i\}.$$
Then $T_1, T_2$ and $T_3$ are three internally edge disjoint $\{x,y,z\}$-trees in $BP_4$.

\vskip 2mm
{\bf Subcase 2.2.}\; For any two integers $a$ and $b$ in $\{j_1, j_2, j_3\}$, it has $a+b\ne 0$.

W.l.o.g., assume that $x\in V(BP_n^1)$, $y\in V(BP_n^{2})$ and $z\in V(BP_n^3)$.

By Lemma \ref{lemBPn1}, there exist vertex subsets $\{u_1,u_2\}\subseteq V(BP_n^1)$, $\{r_1,r_2\}\subseteq V(BP_n^{2})$ and $\{w_1,w_2\}\subseteq V(BP_n^3)$ such that $\{\hat{u}_1, \hat{r}_1, \hat{w}_1\}\subseteq V(BP_n^4)$ and $\{\hat{u}_2, \hat{r}_2, \hat{w}_2\}\subseteq V(BP_n^{\bar{4}})$.

Moreover, there is a vertex $u_3\in V(BP_n^1)\backslash \{u_1,u_2\}$ and a vertex $w_3\in V(BP_n^3)\backslash \{w_1,w_2\}$ such that $\{\hat{u}_3, \hat{w}_3\}\subseteq V(BP_n^{\bar{2}})$. There is a vertex $r_3\in V(BP_n^2)\backslash \{r_1,r_2\}$ such that $\hat{r}_3\in V(BP_n^{\bar{3}})$. Since $BP_n^{\bar{2}}\bigoplus BP_n^{\bar{3}}$ is connected, there is an $\{\hat{u}_3, \hat{w}_3, \hat{r}_3\}$-tree $\widetilde{T}_3$ in $BP_n^{\bar{2}}\bigoplus BP_n^{\bar{3}}$.

For $1\le i\le 3$, let $P_i, Q_i, R_i$ and $T_i$ be the same as in Subcase 2.1. Then $T_1, T_2$ and $T_3$ are three internally edge disjoint $\{x,y,z\}$-trees in $BP_4$.

\vskip 2mm
{\bf Case 3.}\; $n=3$.

Let $H=BP_n^{j_1}\oplus BP_n^{j_2}\oplus BP_n^{j_3}$, where $\{j_1,j_2,j_3\}\subseteq [[n]]$. By Lemma \ref{lemBPn1}, the following Claim \ref{claimBP1} holds clearly.

\begin{claim}\label{claimBP1}
Both $H$ and $BP_3\setminus V(H)$ are connected.
\end{claim}

Recall that $\hat{x}, \hat{y}$ and $\hat{z}$ are out-neighbours of $x,y$ and $z$, respectively. We consider the following cases.

\vskip 2mm
{\bf Subcase 3.1.}\; $\{\hat{x}, \hat{y}, \hat{z}\}\cap V(H)=\emptyset$.

There is an $\{x,y,z\}$-tree $T_1$ in $H$ and an $\{\hat{x}, \hat{y}, \hat{z}\}$-tree $\widetilde{T}_2$ in $BP_3\setminus V(H)$ by Claim \ref{claimBP1}. Let $T_2=\widetilde{T}_2\cup \{x\hat{x}, y\hat{y}, z\hat{z}\}$. Then $T_1$ and $T_2$ are two internally edge disjoint $\{x,y,z\}$-trees in $BP_3$.

\vskip 2mm
{\bf Subcase 3.2.}\; $|\{\hat{x}, \hat{y}, \hat{z}\}\cap V(H)|=1$.

W.l.o.g., assume that $\{\hat{x}, \hat{y}, \hat{z}\}\cap V(H)=\hat{x}$. By Lemma \ref{lemBPn3}, there exists a vertex $u_1\in N_{BP_n^{j_1}}[x]$ satisfying $\hat{u}_1\notin V(H)$.

By Figure \ref{figBP3} and by Lemma \ref{lemBPn1}, the following claim is not difficult to obtained.

\begin{claim}\label{claimBP2}
Let $a_1, a_2$ and $a_3$ be any three vertices such that $a_i\in V(BP_n^{j_i})$ for $1\le i\le 3$. Then $H\setminus \{a_1, a_2,a_3\}$ is connected.
\end{claim}

There is an $\{x,y,z\}$-tree $T_1$ in $H\setminus \{u_1\}$ by Claim \ref{claimBP2} and a $\{\hat{u}_1,\hat{y},\hat{z}\}$-tree $\widetilde{T}_2$ in $BP_3\setminus V(H)$. Let $T_2=\widetilde{T}_2\cup \{xu_1, u_1\hat{u}_1, y\hat{y}, z\hat{z}\}$. Then $T_1$ and $T_2$ are two internally edge disjoint $\{x,y,z\}$-trees in $BP_3$.

\vskip 2mm
{\bf Subcase 3.3.}\; $|\{\hat{x}, \hat{y}, \hat{z}\}\cap V(H)|\ge 2$.

By similar arguments in Subcase 3.2, we can obtain two internally edge disjoint $\{x,y,z\}$-trees in $BP_3$ and the proof is omitted. \hfill$\Box$

\begin{theorem}\label{thmBPn}
For $n\ge 2$, $\kappa_3(BP_n)=n-1$.
\end{theorem}

\noindent{\bf Proof}\; By Lemma \ref{lemupperK} and Lemma \ref{lemBPn1}, we have $\kappa_3(BP_n)\le \delta(BP_n)-1=n-1$. We shall prove the reverse inequality by induction on $n$. It is easily seen from Lemma \ref{lemkr} that $\kappa_3(BP_2)\ge 1$ since $\kappa(BP_2)=2$.  Now suppose that $n\ge 3$ and the theorem holds for any integer $l<n$, i.e., $\kappa_3(BP_l)=l-1$. Let $S=\{x,y,z\}$ be any 3-subset of $V(BP_n)$. The following cases are distinguished.

\vskip 2mm
{\bf Case 1.}\; $x, y$ and $z$ belong to a same cluster of $BP_n$.

W.l.o.g., assume that $\{x, y, z\}\subseteq V(BP_n^1)$. By induction hypothesis, there exist ($n-2$)-internally edge disjoint $S$-trees $T_1, \cdots, T_{n-2}$ in $BP_n^1$ since $BP_n^1$ is isomorphic to $BP_{n-1}$. Recall that $\hat{x}, \hat{y}$ and $\hat{z}$ are out-neighbours of $x, y$ and $z$, respectively. There is a $\{\hat{x}, \hat{y},\hat{z}\}$-tree $\widetilde{T}_{n-1}$ in $BP_n\setminus V(BP_n^1)$ since $BP_n\setminus V(BP_n^1)$ is connected. Let $T_{n-1}=\widetilde{T}_{n-1}\cup \{x\hat{x}, y\hat{y}, z\hat{z}\}$. Then $T_1, \cdots, T_{n-2}, T_{n-1}$ are ($n-1$) desired $S$-trees in $BP_n$.

\vskip 2mm
{\bf Case 2.}\; $x, y$ and $z$ belong to two different clusters of $BP_n$.

W.l.o.g., assume that $\{x, y\}\subseteq V(BP_n^1)$ and $z\in V(BP_n^j)$, where $j\in [[n]]\setminus\{1\}$. There are ($n-1$) internally disjoint ($x,y$)-paths $P_1, \cdots, P_{n-1}$ in $BP_n^1$ by Lemma \ref{lemxypath} and Lemma \ref{lemBPn1}. Let $u_1, \cdots, u_{n-1}$ be neighbour vertices of $x$ such that $u_i\in V(P_i)$ for $1\le i\le n-1$. By Lemma \ref{lemBPn2} and Lemma \ref{lemKfan}, there exists an ($n-1$)-fan $Q_1, \cdots, Q_{n-1}$ in $BP_n\setminus V(BP_n^1)$ from $z$ to $\{\hat{u}_1, \cdots, \hat{u}_{n-1}\}$ such that $\hat{u}_i\in V(Q_i)$ for $1\le i\le n-1$.

For $1\le i\le n-1$, let $T_i=P_i\cup Q_i \cup \{u_i\hat{u}_i\}$. Then $T_1, \cdots, T_{n-1}$ are ($n-1$) desired $S$-trees in $BP_n$.

\vskip 2mm
{\bf Case 3.}\; $x, y$ and $z$ belong to three different clusters of $BP_n$.

According to Lemma \ref{lemBPn4}, there are ($n-1$)-internally edge disjoint $S$-trees in $BP_n$. The proof is completed. \hfill $\Box$

\section{The generalized 3-connectivity of the godan graph}\label{secEAn}

For convenience, we denote the permutation $\scriptsize{\begin{pmatrix} 1&2&\cdots & n\\ p_1&p_2&\cdots & p_n \end{pmatrix}}$ by $p_1p_2\cdots p_n$ and the permutation $\scriptsize{\begin{pmatrix} 1&\cdots& i&\cdots &j \cdots& n\\ 1&\cdots& j&\cdots &i \cdots& n\end{pmatrix}}$ by $(ij)$, while the latter is called a transposition. The composition $\sigma\tau$ of two permutations $\sigma$ and $\tau$ is the function that maps any element $i$ to $\sigma(\tau(i))$.

\begin{figure}[htbp]
\begin{minipage}[t]{0.4\linewidth}
\centering
\resizebox{0.7\textwidth}{!} {\includegraphics{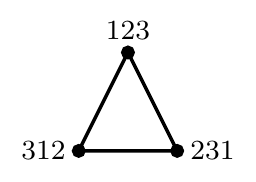}}
\caption{\small The alternating group network $AN_3$} \label{figAN3}
\end{minipage}
%\quad
\begin{minipage}[t]{0.6\linewidth}
\centering
\resizebox{0.7\textwidth}{!} {\includegraphics{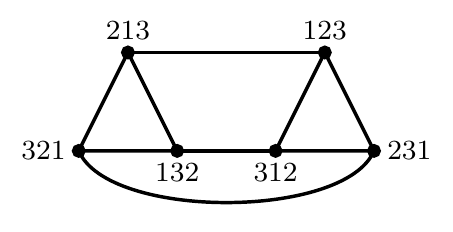}}
\caption{\small The godan graph $EA_3$} \label{figEA3}
\end{minipage}
\end{figure}

Recall that $S_n$ is the symmetric group on $[n]$. The alternating group $A_n$ ($n\ge 3$) is the subgroup of $S_n$ containing all even permutations. Let $\Omega=\{(123),(132)\}\cup \{(12)(3i)| 4\le i\le n\}$. It is well known that $\Omega$ is a generating set for $A_n$. For $n\ge 3$, the $n$-dimensional alternating group network $AN_n$ (\cite{Ji1998}) is a graph with vertex set $V(AN_n)=A_n$ and edge set $E(AN_n)=\{(u,v)| u=vs, \;{\rm where}\; \{u,v\}\subseteq V(AN_n) \;{\rm and}\; s\in \Omega\}$. It is seen that the $n$-dimensional alternating group network $AN_n$ is the Cayley graph Cay($A_n, \Omega$). The alternating group network $AN_3$ is depicted in Figure \ref{figAN3}.

\begin{lemma}\label{lemANn}(\cite{Ji1998,ANn12018,ANn12019,ANn2018})
For the $n$-dimensional alternating group network $AN_n$ ($n\ge 3$), the following properties hold:\\
(1) $AN_n$ is ($n-1$)-regular with $\frac{n!}{2}$ vertices and $\frac{n!(n-1)}{4}$ edges;\\
(2) $\kappa(AN_n)=n-1$;\\
(3) $\kappa_3(AN_n)=n-2$.
\end{lemma}

\begin{definition}\label{defEAn}(\cite{GodanIn2022,godangraph2022})
For $n\ge 3$, let $\Omega^*=\{(12),(123),(132)\}\cup \{(12)(3i)|4\le i\le n\}$. The $n$-dimensional godan graph $EA_n$ is the graph with vertex set $V(EA_n)=S_n$ in which two vertices $u$ and $v$ are adjacent in $EA_n$ if and only if $u=vs^*$, where $s^*\in \Omega^*$.
\end{definition}
It is easily seen that $EA_n$ is isomorphic to the Cayley graph Cay($S_n, \Omega^*$). The godan graph $EA_3$ is depicted in Figure \ref{figEA3}.

\begin{lemma}\label{lemEAn}(\cite{GodanIn2022,godangraph2022})
For $n\ge 3$, the $n$-dimensional godan graph $EA_n$ has the following properties:\\
(1) it is $n$-regular with $n!$ vertices and $\frac{n\times n!}{2}$ edges;\\
(2) $EA_n[A_n]\cong AN_n$, $EA_n[S_n\backslash A_n]\cong AN_n$;\\
(3) $E(A_n, S_n\backslash A_n)$ is a perfect matching of $EA_n$, where $E(A_n, S_n\backslash A_n)$ denotes all cross edges between $A_n$ and $S_n\backslash A_n$;\\
(4) $\kappa(EA_n)=n$.
\end{lemma}

For simplicity, we denote $EA_n[A_n]$ by $AN_n^1$ and denote $EA_n[S_n\backslash A_n]$ by $AN_n^2$, respectively. Note that $AN_n^i$ is isomorphic to $AN_n$ for $1\le i\le 2$. Then we can write the construction of $EA_n$ symbolically as
$$EA_n=AN_n^1\oplus AN_n^2.$$
For $1\le i\le 2$, let $x$ be a vertex in $AN_n^i$, it is clearly that $x$ has an out-neighbor $\hat{x}\in V(AN_n^{3-i})$. Now we shall investigate the generalized 3-connectivity of $EA_n$.

\begin{theorem}\label{thmEAn}
For $n\ge 3$, $\kappa_3(EA_n)=n-1$.
\end{theorem}

\noindent{\bf Proof}\; Since $EA_n$ is $n$-regular, we have that $\kappa_3(EA_n)\le n-1$ by Lemma \ref{lemupperK}. To prove the theorem, it suffices to show that there are ($n-1$)-internally edge disjoint $S$-trees for any $S=\{x,y,z\}\subseteq V(EA_n)$. According to the construction of $EA_n$, the following cases are considered.

\vskip 2mm
 {\bf Case 1.} $|S\cap V(AN_n^1)|=3$.

That means $\{x,y,z\}\subseteq V(AN_n^1)$. It follows from Lemma \ref{lemANn} and Lemma \ref{lemEAn} that there are ($n-2$)-internally edge disjoint $S$-trees $T_1, \cdots, T_{n-2}$ in $AN_n^1$ since $\kappa_3(AN_n^1)=\kappa_3(AN_n)=n-2$.

By Lemma \ref{lemANn} and Lemma \ref{lemEAn}, there exists an $\{\hat{x}, \hat{y}, \hat{z}\}$-tree $\widetilde{T}_{n-1}$ in $AN_n^2$ since $AN_n^2$ is connected. Let $T_{n-1}=\widetilde{T}_{n-1}\cup\{x\hat{x}, y\hat{y}, z\hat{z}\}$. It is clearly that $T_1, \cdots, T_{n-2}, T_{n-1}$ are ($n-1$)-internally edge disjoint $S$-trees in $EA_n$.

\vskip 2mm
{\bf Case 2.} $|S\cap V(AN_n^1)|=2$.

W.l.o.g., assume that $S\cap V(AN_n^1)=\{x, y\}$ and $S\cap V(AN_n^2)=\{z\}$. By Lemma \ref{lemANn}, Lemma \ref{lemEAn} and Menger's theorem \cite{Bondy}, there are ($n-1$) internally disjoint $(x,y)$-paths $P_1, \cdots, P_{n-1}$ in $AN_n^1$.  For $1\le i\le n-1$, let $w_i$ be a neighbour vertex of $x$ in $P_i$, moreover, let $\hat{w}_i$ be the out-neighbour of $w_i$.

Let $\hat{W}=\{\hat{w}_1, \cdots, \hat{w}_{n-1}\}$.

By Lemma \ref{lemKfan}, there are ($n-1$) internally disjoint ($z, \hat{W}$)-paths $\widetilde{P}_1, \cdots, \widetilde{P}_{n-1}$ in $AN_n^2$ where $\hat{u}_i\in V(\widetilde{P}_i)$ for $1\le i\le n-1$. It is possible that $z\in \hat{W}$. We can let $\widetilde{P}_i=\{z\}$ if $z=\hat{w}_i$ for $1\le i\le n-1$.

For $1\le i\le n-1$, let $T_i=P_i\cup\widetilde{P}_i \cup \{w_i\hat{w}_i\}$. Then $T_1, \cdots, T_{n-2}, T_{n-1}$ are ($n-1$)-internally edge disjoint $S$-trees in $EA_n$.

%If $z\in \hat{W}$. W.l.o.g., assume that $z=\hat{w}_{n-1}$. By Lemma \ref{lemKfan}, there are ($n-2$) internally disjoint ($z, \hat{W}-\{z\}$)-paths $\widetilde{P}_1, \cdots, \widetilde{P}_{n-2}$ in $AN_n^2$. For $1\le i\le n-2$, let $T_i=P_i\cup\widetilde{P}_i \cup \{w_i\hat{w}_i\}$. Let $T_{n-1}=P_{n-1}\cup\{w_{n-1}z\}$.
%

\vskip 2mm
{\bf Case 3.} $|S\cap V(AN_n^1)|=1$.

That means $|S\cap V(AN_n^2)|=2$. By Lemma \ref{lemEAn}, $AN_n^2$ is isomorphic to $AN_n$. We can get that there are ($n-1$)-internally edge disjoint $S$-trees in $EA_{n}$ by similar arguments in Case 2.

\vskip 2mm
 {\bf Case 4}. $|S\cap V(AN_n^1)|=0$.

That means $\{x,y,z\}\subseteq V(AN_n^2)$. By similar arguments in Case 1, we conclude that there are ($n-1$)-internally edge disjoint $S$-trees in $EA_{n}$.  \hfill$\Box$

\section{Conclusion}

The generalized $k$-connectivity is a natural generalization of the classical connectivity and can serve for measuring the capability of a network $G$ to connect any $k$ vertices in $G$. In this paper, we focus on the generalized 3-connectivity of the burnt pancake graph $BP_n$ and the $n$-dimensional godan graph $EA_n$.  We can see that most of the results on the generalized $k$-connectivity of networks are about $k=3$. It would be an interesting and challenging topic to study $\kappa_k(BP_n)$ and $\kappa_k(EA_n)$ for $k\ge 4$.

%\begin{acknowledgements}
%
%\end{acknowledgements}



\end{document}